\documentclass[11pt]{article}
\usepackage{amssymb,amsmath,amsthm}
\usepackage{epsfig}
\usepackage{psfrag}
\newtheorem{Theorem}{\sc Theorem}[section]
\newtheorem{Lemma}[Theorem]{\sc Lemma}
\newtheorem{Proposition}[Theorem]{\sc Proposition}
\newtheorem{Corollary}[Theorem]{\sc Corollary}
\newtheorem{Definition}[Theorem]{\sc Definition}
\newtheorem{Example}[Theorem]{\sc Example}
\newtheorem{Remark}[Theorem]{\sc Remark}

\def\hpic #1 #2 {\mbox{$\begin{array}[c]{l} 
\epsfig{file=#1,height=#2}\end{array}$}}
\def\wpic #1 #2 {\mbox{$\begin{array}[c]{l} 
\epsfig{file=#1,width=#2}\end{array}$}}

\def\C{\mathbb C}

\def\R{\mathbb R}

\def\CE{{\cal {E}}}
\def\CF{{\cal {F}}}

\def\CH{{\cal {H}}}

\def\CK{{\cal {K}}}
\def\CL{{\cal {L}}}

\def\CP{{\cal {P}}}

\def\CS{{\cal {S}}}

\def\CU{{\cal {U}}}

\def\be{\begin{equation}}
\def\ee{\end{equation}}
\def\bt{\begin{Theorem}}
\def\et{\end{Theorem}}
\def\bi{\begin{itemize}}
\def\ei{\end{itemize}}
\def\bea{\begin{eqnarray}}
\def\eea{\end{eqnarray}}
\def\ba{\begin{array}}
\def\ea{\end{array}}
\def\beast{\begin{eqnarray*}}
\def\eeast{\end{eqnarray*}}
\def\ben{\begin{enumerate}}
\def\een{\end{enumerate}}

\def\bi{\bibitem}

\newcommand{\e }{\mbox{$ \epsilon $}}

\def\rar{\rightarrow}
\def\Rar{\Rightarrow}

\def\Lra{{\Leftrightarrow}}

\def\e{{\epsilon}}

\begin{document}
\begin{center}
{\Large {\bf Hilbert von Neumann  modules}}
\end{center}

\bigskip

\begin{center}
Panchugopal Bikram\footnote{The Institute of Mathematical Sciences, Chennai},
Kunal Mukherjee$^1$, R. Srinivasan\footnote{Chennai Mathematics                 
  Institute, Chennai},\\
and V.S. Sunder$^1$
\end{center}

\begin{abstract}
We introduce a way of regarding Hilbert von Neumann modules as spaces
of operators between Hilbert space, not unlike [Skei], but in an
apparently much simpler manner and involving far less machinery. We verify 
that our definition is equivalent to that of [Skei], 
by verifying the `Riesz lemma' or what is called `self-duality' in
[Skei]. An advantage with our approach is that we can totally
side-step the need to go through $C^*$-modules and avoid the two
stages of completion - first in norm, then in the strong operator
topology - involved in the former approach.  

We establish the analogue of the Stinespring dilation theorem for
Hilbert von Neumann bimodules, and we develop our version of `internal
tensor products' which we refer to as Connes fusion for obvious
reasons.

In our discussion of examples, we examine the bimodules arising from
automorphisms of von Neumann algebras, verify that fusion of bimodules
corresponds to composition of automorphisms in this case, and that the
isomorphism class of such a bimodule depends only on the inner conjugacy
class of the automorphism. We also
relate Jones' basic construction to the Stinespring dilation
associated to the conditional expectation onto a finite-index
inclusion (by invoking the uniqueness assertion regarding the latter).

\end{abstract}

\section{Preliminaries}

The symbols $\CH$ and $\CK$, possibly anointed with subscripts or
other decorations, will always denote complex separable Hilbert
spaces, while $\CL(\CH,\CK)$ will denote the set of bounded operators
from $\CH$ to $\CK$. For $E \subset \CL(\CH,\CK)$, we shall write
$\left[E\right]$ for the closure, in the weak operator topology (WOT,
in the sequel), of the linear subspace of $\CL(\CH,\CK)$ spanned by
$E$. Similarly, if $\CS \subset \CH$ is a set of vectors, we shall
write $\left[\CS\right]$ for the norm-closed subspace of $\CH$ spanned
by $\CS$. 

Without explicitly citing it again to justify statements we make, we
shall use the fact that a linear subspace of $\CH$ (resp.,
$\CL(\CH,\CK)$) is closed in the weak topoogy (resp., WOT) if and only
if it is closed in the strong or norm topology (resp., `SOT'). (For
example, $[E]$ is an algebra if $E$ is.)

If $E \subset \CL(\CH,\CK)$ and $F \subset \CL(\CH_1,\CH)$, we write 
\[ EF = \{xy : x \in E, y \in F\} \mbox{ and } E^* = \{x^*: x \in E\}~.\]

If $i:\CH_0 \hookrightarrow \CH$ and $j:\CK_0 \hookrightarrow \CK$, then
we shall think of $\CL(\CH,\CK_0)$ as the subset $f \CL(\CH_0,\CK) e = 
j \CL(\CH_0,\CK) i$ of 
$\CL(\CH_0,\CK)$, where $e$ and $f$ are the projections $e=i^*, f=j^*$.

\begin{Proposition}\label{vncorner}
For $i=1,2$, let $e_i$ denote the projection of  $\CH_1 \oplus \CH_2$
onto $\CH_i$. The following conditions on an $E \subset \CL(\CH_2,\CH_1)$
are equivalent:
\ben
\item There exists a von Neumann algebra $M \subset \CL(\CH_1 \oplus
  \CH_2)$ such that $e_1,e_2 \in M$ and $E = e_1Me_2$.
\item $E = \left[E\right] \supset EE^*E.$
\een 

When these equivalent conditions are met, we shall say that
$(E,\CH_1,\CH_2)$ is a {\bf (1,2) von Neumann corner}.
\end{Proposition}

\begin{proof}
$(1) \Rar (2)$  is obvious.

$(2) \Rar (1)$: Observe that the assumption (2) implies that $\left[E^*E\right]$
  is a WOT-closed *-subalgebra of $\CL(\CH_2)$. Let $p_2 
  = \sup \{ p: p \in \CP(\left[E^*E\right]\}$ and define $M_{22} =
  [E^*E] + \C(e_2 - p_2)$; so $M_{22}$ is a von Neumann subalgebra of
  $\CL(\CH_2)$ and $e_2 - p_2$ is a central minimal projection in it.

Similarly,  define $M_{11} = [EE^*] + \C(e_1 - p_1)$, where $p_1 
  = \sup \{ p: p \in \CP(\left[EE^*\right]\}$; so $M_{11}$ is a von
  Neumann subalgebra of $\CL(\CH_1)$  and $e_1 - p_1$ is a central
  minimal projection in it.

Finally set $M_{12} = E, M_{21} = E^*$ and $M = \sum_{i,j=1}^2 M_{ij}$.  (Alternatively $M$ is the von Neumann algebra $(E \cup E^*)^{\prime \prime}$ ; and it is clear that $E =  e_1Me_2$.

\end{proof}

\begin{Definition}\label{nondegdef}
{\rm 
\ben
\item The projection $p_1$ (resp. $p_2$) ocurring in the proof of Proposition \ref{vncorner} will be referred to as the {\bf left-support} (resp., {\bf right-support}) {\bf projection} of the (1,2) von Neumann corner $E$.
\item A (1,2) von Neumann corner $(E,\CH_1,\CH_2)$ will be said to be
 {\bf non-degenerate} if its support projections are as large as they
 can be:  i.e., $p_i = (e_i =) 1_{\CH_i}, i=1,2$.
 \een }
\end{Definition}

\begin{Remark}\label{nondegrmk}
\ben
\item {\rm The support projections $p_1, p_2$ of $E$ have the following equivalent descriptions:
\begin{itemize}
\item $ran~p_1 = \left[\bigcup\{ran ~x: x \in E\}\right] = \left( \bigcap \{ker~x^*: x \in E\}^\perp \right)$; and
\item $ran~p_2 = \left[\bigcup\{ran ~x^*: x \in E\}\right] = \left( \bigcap \{ker~x: x \in E\}^\perp \right)$.

\end{itemize}
\item A (1,2) von Neumann corner $(E,\CH_1,\CH_2)$ is non-degenerate precisely when
  $M_{11} (E) = \left[EE^*\right]$ and $M_{22}(E) =
  \left[E^*E\right]$ are unital von Neumann subalgebras of
 $\CL(\CH_1)$   and $\CL(\CH_2)$ respectively. }
\een 

\end{Remark}

\begin{Definition}\label{HvNmoddef}
\ben
\item  If $A_2$ is a von Nemann algebra, a {\bf Hilbert von Neumann 
  $A_2$ - module} is a tuple $\CE = (E,\CH_1,(\pi_2,\CH_2))$ where
  $(E,\CH_1,\CH_2)$ is a (1,2) von
  Neumann corner equipped with a normal {\rm isomorphism} $\pi_2:A_2
  \rar \left[E^*E\right]$.  
\item A {\bf submodule} of a Hilbert von Neumann 
  $A_2$-module $E$ is a  subset $E_1 \subset E$ satisfying
  \[E_1 = \left[E_1\right] \supset E_1E^*E.\]
\item If $A_1, A_2$ are von Neumann algebras, a {\bf Hilbert von Neumann 
  $A_1-A_2$ - bimodule} is a tuple
\[\CE = (E,(\pi_1,\CH_1),(\pi_2,\CH_2))\]
comprising a Hilbert von Neumann 
  $A_2$ - module $(E,\CH_1,(\pi_2,\CH_2))$ equipped with a
  normal unital homomorphism $\pi_1:A_1 \rar \left[EE^*\right]$ (where the `unital requirement' is that $\pi_1(1_{A_1}) = p_1$ is the identity of $[EE^*]$).
\een
\end{Definition}
 
\begin{Remark}\label{HvNmodrmks}
{\rm \ben
\item If $E \subset \CL(\CH_2, \CH_1)$ is any (possibly degenerate)
  (1,2) von Neumann corner, 
with associated support projections  $p_1, p_2$ (as in  Definition
\ref{nondegdef}), define $\CK_i = ran ~p_i, A_1 =  [EE^*], A_2 = [E^*E]$
and let $\pi_i$ denote the identity representation of $A_i$ on
$\CK_i$; ithen  $(E, (\pi_1, \CK_1), (\pi_2, \CK_2))$ is seen to be a
non-degenerate Hilbert von Neumann $A_1-A_2$ - bimodule. This is why
non-degeneracy is not a serious restriction. 
\item
A Hilbert von Neumann  $A_2$ - module $(E,\CH_1,(\pi_2,\CH_2))$
does indeed admit a right-$A_2$ action and an
$A_2$ - valued inner product thus:
\[x \cdot a_2 = x\pi_2(a_2) ~; \langle x_1, x_2 \rangle_{A_2}
= \pi_2^{-1}(x_1^*x_2) \]
(Here and in the sequel, we shall write $\langle \cdot,\cdot \rangle_B$ for the $B$ - valued inner-product on a Hilbert $B$ - module.)
Notice, further, that the norm $E$ acquires from this Hilbert $A_2$ - module
structure is nothing but the operator norm on $E$.
\item A submodule of a Hilbert von Neumann  $A_2$ - module is a (possibly degenarate) (1,2) von Neumann corner.
\item
In a general  Hilbert von Neumann  $A_2$ - module $\CE = (E,\CH_1,(\pi_2,\CH_2))$, note that 
\[[EE^*]  \ni a \mapsto \left( E \ni x \mapsto a \cdot x =: ax  \right)\]
defines a *-homomorphism  of $[EE^*]$ into the space $\CL^a(E)$ of bounded
adjointable operators on $E$, since, for instance
\beast \langle a \cdot x, y \rangle_{A_2} &=& (ax)^*y\\
&=& x^* (a^*y)\\
&=& x^* (a^* \cdot y)\\
&=& \langle x, a^* \cdot y \rangle_{A_2}~.
\eeast
\item In the language of (2) above, the `rank-one operator'
  $\theta_{x,y}$ is seen to be given by \beast
  \theta_{x,y}(z) &=& x \langle y,z\rangle_{A_2}\\
  &=& xy^*z~, \eeast so that the `rank-one operator' $\theta_{x,y}$ on
  $E$ is nothing but left multiplication by $xy^*$ on $E$, for any
  $x,y \in E$. Let us write $B = [EE^*]$, $C = A_2$ and $A$ for the
  norm-closure of the linear span of $EE^*$.  Then it is clear that
  $A$ is a norm-closed ideal in $B$, and that there is a unique
  $C^*$ - algebra isomorphism $\alpha:A \rar \CK(E)$ such that
  $\alpha(xy^*) = \theta_{x,y}, ~\forall x,y \in E.$ If $E$ is
  non-degenerate, then $A$ is an essential ideal in $B$ and $\alpha$
  is injective. It then follows from [Lan] Proposition 2.1, that
  $\alpha$ extends uniquely to an isomorphism of $B$ onto $\CL^a(E)$.
  (In fact, the reason for introducing the symbols $A,B,C$ above was
  in order to use exactly the same symbols as in the Proposition 2.1
  referred to above.)  
  \item This remark concerns our requirement, in the definition of a Hilbert von Neumann $A_2$-module, that $\pi_2:A_2 \rar [E^*E]$ must be an isomorphism. What is really needed is that $\pi_2$ is onto. If $\pi_2$ is merely surjective but not injective, there must exist a central projection $z \in A_2$ such that $ker~\pi_2 = (1-z) A_2$ so $\pi_2$ would map $zA_2$ isomorphically onto $[E^*E]$ and the $A_2$-valued inner product (see item (2) of this remark) would actually take values in $zA_2$ and we could apply our analysis to $zA_2$ and think of $A_2$ as acting via its quotient (and ideal) $zA_2$.
\item The `unital requirement' made in the definition of a Hilbert von Neumann bimodule has the consequence that $\pi_1(A_1)E = E$.
  \een }
\end{Remark}

\begin{Lemma}\label{Epd}
Let $(E,\CH_1,\CH_2)$ be a (1,2) von Neumann corner. Suppose $x \in
\CL(\CH_2,\CH_1)$ has polar decomposition $x = u|x|$. Then the
following conditions are equivalent: 
\ben
\item $x \in E$.
\item $u \in E \mbox{ and } |x| \in [E^*E]$.
\item $u \in E \mbox{ and } |x^*| \in [EE^*]$.
\een
\end{Lemma}

\begin{proof} Since $(2) \Rar (1)$ and $(3) \Rar (1)$ are obvious, let
  us prove the reverse implications. So,
  suppose $x \in E$. Then $x^*x \in E^*E$ (resp., $xx^* \in EE^*$) and
  as, $|t|$ is uniformly approximable on compact subsets of $\R$ by 
  polynomials with vanishing constant term, it is seen that $|x| \in
  \left[E^*E\right]$ and $|x^*| \in \left[EE^*\right]$.
Define $f_n \in C_0([0,\infty))$ by
\[f_n(t) = \left\{ \ba{ll} 0 &  \mbox{if } t < \frac{1}{2n}\\
2n^2(t- \frac{1}{2n}) & \mbox{if } \frac{1}{2n} \leq t \leq \frac{1}{n}\\
\frac{1}{t} & \mbox{if } t \geq \frac{1}{n} \ea \right.\]
Since $f_n$ is uniformly approximable on $sp(|x|)$ by
polynomials with 
vanishing constant term, it is seen that $f_n(|x|) \in
\left[E^*E\right]$, and 
hence $xf_n(|x|) \in E$. It follows from
the definitions that 
$|x|f_n(|x|)$ WOT-converges to
$1_{(0,\infty)}(|x|) = u^*u $. In particular, 
$u = u(u^*u) = WOT-lim ~u(|x|f_n(|x|) = WOT-lim ~xf_n(|x|) \in
[E\left[E^*E\right]] \subset E$.

\end{proof}

\begin{Proposition}\label{submodperp} 
If $E_1$ is a submodule of a Hilbert von Neumann $A_2$ - module $E$, and if $E_1 \neq E$, there exists a non-zero $y \in E$ such that $y^*x = 0 ~\forall x \in E_1$.
\end{Proposition}

\begin{proof}
As observed in Remark \ref{HvNmodrmks}(3), $E_1$ is a possibly
degenerate (1,2) von Neumann corner in  $\CL(\CH_2,\CH_1)$. Let
$p_1=\bigvee \{e : e \in \CP([E_1^*E_1])\}$ and $q_1 =\bigvee \{f : f
\in \CP([E_1E_1^*])\}$ be the right-  and left- support projections of
$E_1$. Similarly, let $p=\bigvee \{e : e \in \CP([E^*E])\}$ and $q
=\bigvee \{f : f \in \CP([EE^*])\}$ be the right-  and left- support
projections of $E$.  

First observe that the hypotheses imply that
\[(E^*E)(E_1^*E_1)(E^*E) = (E_1E^*E)^*(E_1E^*E) \subset E_1^*E_1\]
and hence that $[E_1^*E_1]$ is a WOT-closed ideal in the von Neumann subalgebra $[E^*E]$ of
$\CL(p\CH_2)$;
consequently $p_1=\bigvee \{e : e \in \CP([E_1^*E_1])\}$ is a central projection in $[E^*E]$ and $[E_1^*E_1] = [E^*E]p_1$. It follows that if $x_1 \in E_1$ has polar decomposition $x_1 = u_1|x_1|$, then (by Lemma \ref{Epd}) $u_1 \in E_1$ and $|x_1| \in [E_1^*E_1] = [E^*E]p_1$, and in particular,
$x_1p_1 = u_1|x_1|p_1 =  u_1|x_1| = x_1$; i.e., $E_1 = E_1p_1$.

Next, by definition, $[\bigcup \{ran~x_1: x_1 \in E_1\}] = [\bigcup \{ran~x_1x_1^*: x_1 \in E_1\}] = 
[\bigcup \{ran 1_{(0,\infty)}(|x_1^*|): x_1 \in E_1\}] = ~ran~q_1$; hence if $x_1 \in E_1$, then $x_1 = q_1x_1$, and we see that $E_1 = q_1E_1$.

Summarising the previous two paragraphs, we have
\be \label{e1qp} E_1 = q_1E_1 = E_1p_1 ~.\ee

(In fact, $x_1 = x_1p_1 = q_1x_1 ~\forall ~x_1 \in E_1$.)

We now consider three cases:

\medskip {\em Case 1:} $p_1 \neq p$

Here $(p-p_1) \neq 0$ and the definition of $p$ implies that there exists a $y \in E$ such that $y  = y(p-p_1) \neq 0$. Then, for any $x \in E_1$, we have $x = xp_1$ and hence
\[ y^*x = (p-p_1)y^*x = (p-p_1) y^*x p_1 \in (p-p_1) E^*E p_1 = (p-p_1)p_1E^*E  = \{0\}~.\]

\medskip {\em Case 2:} $q_1 \neq q$

Here $(q-q_1) \neq 0$ and the definition of $q$ implies that there exists a $y \in E$ such that $y  = (q-q_1)y \neq 0$. Then, for any $x \in E_1$, we have $x = q_1x$ and hence
\[ y^*x = y^*(q-q_1)x = y^*(q-q_1)q_1x = 0~.\]

\medskip {\em Case 3:} $p_1 = p, q_1= q$.

We shall show that the hypotheses of this case imply that $E_1 = E$ and hence cannot arise.
To see this, begin by noting that the collection of non-zero partial isometries in $E_1$ is non-empty in view of Lemma \ref{Epd}. (Otherwise $E_1 = \{0\}, p_1=q_1=0$ and so $E = \{0\} = E_1$.) Hence the family $\CF$ of collections $\{u_i: i \in I\}$ of partial isometries in $E_1$ with pairwise orthogonal ranges, is non-empty. Clearly $\CF$ is partially ordered by inclusion, and it is easy to see that Zorn's lemma is applicable to $\CF$. 

If $\{u_i: i \in I\}$ is a maximal element of $\CF$, we assert that $\sum_{i\in I} u_iu_i^* =q$. Indeed, if $(q - \sum_{i\in I} u_iu_i^*) \neq 0$, the assumption $q=q_1$ will imply the existence of an $x_1 \in E_1$ such that $x_1 = (q - \sum_{i\in I} u_iu_i^*)x_1 \neq 0$. Then $x_1 \in [E_1E_1^*E_1] \subset E_1$ and so if $x_1 = v_1|x_1|$ is its polar decomposition, then $v_1 \in E_1 \setminus \{0\}$ and $ran~v_1 = \overline{ran~x_1}$ is orthogonal to $ran~u_i$ for each $i \in I$, thus contradicting the maximality of $\{u_i: i \in I\}$.

Thus, indeed $q = \sum_{i \in I} u_iu_i^*, ~u_i \in E_1$.

Now, if $x \in E$ is arbitrary, then,
\beast
x &=& qx\\
&=& \sum_{i \in I} u_iu_i^*x\\
&\in& [E_1E_1^*E]\\
&\subset& [E_1 E^*E]\\
&\subset& E_1
\eeast
and so $E = E_1$ in this case, and the proof of the Proposition is complete.
\end{proof}

Given a submodule $E_1$ of a Hilbert von Neumann module $E$, as above, we shall write $E_1^\perp$ for the set $\{y \in E : y^*E_1 = \{0\}\}$ and refer to it as the {\bf orthogonal complement of $E_1$ in $E$}. We now reap the consequences of Proposition \ref{submodperp} in the following Corollary.

\begin{Corollary}\label{perpcor}
Let $E_1$ be a submodule of a Hilbert von Neumann $A_2$ - module. Then,
\ben
\item $E_1^\perp = (1-q_1)E$, where $q_1$ is the left support projection of $E_1$.
\item $E_1^{\perp\perp} = q_1E$.
\item If $S$ is any subset of $E$, then $S^{\perp\perp} = [S[E^*E]]$.
\item If $E_1$ is a submodule of a Hilbert von Neumann module $E$, there exists a projection $ q_1 \in [EE^*]$ such that $E_1 = E_1^{\perp \perp} = q_1E$ and $E_1^\perp = (1-q_1)E$; and in particular $E_1$ is complemented in the sense that $E = E_1 \oplus E_1^\perp$.
\een
\end{Corollary}

\begin{proof}
It is clear that $y^*x = 0$ if and only if $y$ and $x$ have mutually orthogonal ranges. 

(1) The previous sentence and the definition of $q_1$ imply that 
\[y \in E_1^\perp \Lra \left(q_1y = 0 \mbox{ and } y \in E \right) \Lra y \in (1-q_1)E.\]

(2) follows from (1) and the definition of the orthogonal complement.

(3) Let $E_1 = [SE^*E]$. It should be clear that $y \in S^\perp \Lra  y \in E_1^\perp = q_1E$,
 by part (1) of this Corollary, and hence that
 \[S^{\perp\perp} = E_1^{\perp\perp} ~.\]
In view of Remark \ref{HvNmodrmks}(1) we may view $S^{\perp \perp}$ as a Hilbert von Neumann bimodule, and regard $E_1$ as a submodule of $S^{\perp \perp}$. We may then deduce from Proposition \ref{submodperp} that if $E_1$ were not equal to $S^{\perp \perp}$,  then there would have to exist a non-zero $y \in S^{\perp \perp}$ such that $y^*E_1 = \{0\}$. This would imply that $y \in S^\perp$ and $y \in S^{\perp\perp}$ so that $y^*y = 0$, a contradiction.

(4) follows from the preceding parts of this Corollary.
\end{proof}

That our definitions of Hilbert von Neumann  modules and bimodules are
consistent with those of [Skei] is a
consequence of the following 
version of Riesz' Lemma, which establishes that our Hilbert von
Neumann  modules are indeed `self-dual' which is one of the equivalent
conditions for a von Neumann module in the sense of [Skei].

On the other hand, it is clear from [Skei] that any Hilbert von Neumann
$A_2$ - module in the sense of [Skei] is also a Hilbert von Neumann
$A_2$ - module in our sense, and the two formulations are thus equivalent.

\begin{Proposition}\label{rielem}({\bf Riesz lemma})
Suppose $\CE$ is a Hilbert von Neumann $A_2$ - module, and $f: E
\rar A_2$ is right $A_2$-linear - meaning $f(x\pi_2(a_2)) =
\pi_2^{-1}(f(x)\pi_2(a_2))$ for all $x \in E, a_2 \in A_2$, or equivalently and less clumsily,  suppose $f:E \rar [E^*E]$ is linear and satisfies
$f(xz) = f(x)z$ for all $x \in E, ~z \in [E^*E]$; and suppose $f$ is bounded -
meaning $\|f(x)\| \leq K\|x\|$ for all $x \in E$, and some $K > 0$. Then there exists
$y \in E$ such that $f(x) =  y^*x ~\forall x \in E$.
\end{Proposition}

\begin{proof}
First notice that if $x \in E$ has polar decomposition $x=u|x|$ (so $u
\in E, |x| \in [E^*E] = \pi_2(A_2)$, and if $\xi \in \CH_2$, then
\bea\label{tfineq}
\|f(x)\xi\| &=& \|f(u)|x|\xi\| ~~~\mbox{(by right $A_2$ - linearity of
  $f$)}\nonumber\\
&\leq& \|f(u)\| \||x|\xi\|\nonumber\\
&\leq& K \||x|\xi\|\nonumber\\
&=& K \|u^*x\xi\|\nonumber\\
&\leq& K \|x\xi\| ~.
\eea

Next, find vectors $\xi_n \in \CH_2$ such that $\CH_2 = \oplus_n
[\pi_2(A_2)\xi_n]$ (orthogonal direct sum). It follows that $p_1\CH_1 = \oplus_n [E\xi_n]$, where $p_1$ is the left-support projection of $E$ (because
if $n \neq m$ and $x,y \in E$, then
\[ \langle x\xi_n,y\xi_m \rangle = \langle \xi_n, x^*y\xi_m \rangle =
0\]
and \[ [\bigcup_n [E\xi_n]] = [\bigcup_n [EE^*E\xi_n]] = 
[E\CH_2] = p_1\CH_1~.\]

Infer from the above paragraph and equation \ref{tfineq} that for
arbitrary $a_n \in A_2$ with $\sum_n \|\pi _2(a_n)\xi_n)\|^2 < \infty$ and $x \in E$, we have
\beast
\|f(x)(\sum_n \pi_2(a_n)\xi_n)\|^2 &=& \|\sum_n (f(x)\pi_2(a_n))\xi_n\|^2\\
&=& \sum_n \|f(x\pi_2(a_n))\xi_n\|^2\\
&\leq& \sum_n K^2 \|x\pi_2(a_n)\xi_n\|^2  ~~~~\mbox{(by eq. (\ref{tfineq})}\\
&=& K^2 \|x(\sum_n \pi_2(a_n)\xi_n)\|^2 ~;
\eeast

Now deduce that there exists a
unique bounded operator $z_f \in \CL(\CH_1,\CH_2)$ satisfying $z_f = z_fp_1$ and
\[z_f(x\xi) = f(x)\xi ~, \forall x \in E, \xi \in \CH_2~.\]

The definition of $z_f$ implies that $z_fE \subset [E^*E]$; hence
\[z_f = z_f p_1 \in z_f[EE^*] \subset [z_f EE^*] \subset [[E^*E]E^*] = E^*~.\]

So $y =: z_f^* \in E$ and we have
\[f(x) = z_fx = y^*x \]
as desired.
\end{proof}

\section{Standard bimodules and complete positivity}
Given an element $x$ of a von Neumann algebra $M$, et us write $pr(x)$
for the projection onto the range of $x$. (Thus $pr(x) = 1_{(0,\infty)}(xx^*)$.)

\begin{Lemma}\label{cp}
Suppose $\eta :A \rar B$ is a normal positive linear map of von
Neumann algebras. Let $e_\eta = \bigvee\{u~pr(\eta(1))~u^*:
u \in \CU(B)\}$ be the ($B$-)central support of $pr(\eta(1)$. Then
the smallest WOT-closed ideal in $B$ which 
contains $\eta(A)$ (equivalently $\eta(1)$) is $e_\eta B$. (In particular, $\eta(a) = e_\eta
\eta(a) ~\forall a \in A$.)
\end{Lemma}

\begin{proof} 
If $p \in \CP(A)$, then $\eta(p) \leq \eta(1) ~\Rar pr(\eta(p) \leq
pr(\eta(1) \leq e_\eta$. Hence $\eta(p) = e_\eta \eta(p) \in e_\eta B$, so also
$B\eta(p)B \subset e_\eta B$. Conclude that $[B\eta(A)B] = [B\eta([\CP(A)])B] =
[B\eta(\CP(A))B] \subset e_\eta B$. Conversely,
$[B\eta(A)B] \supset [B\CU(B)\eta(1)\CU(B)B] \supset [B e_\eta B] =
e_\eta B$, and the proof is complete.
\end{proof}

\begin{Definition}
A Hilbert von Neumann $A_2$ - module $\CE = (E,\CH_1,(\pi_2,\CH_2))$ will be called {\bf standard} if :
\begin{itemize}
\item $\CH_2 = L^2(A_2,\phi)$ for some faithful normal state $\phi$ on $A_2$;
\item $\pi_2$ is the left-regular representation; and
\item $E$ is non-degenerate.
\end{itemize}

A Hilbert von Neumann $A_1-A_2$ - bimodule will be called standard if it is standard as a Hilbert von Neumann $A_2$ - module.
\end{Definition}

\bt \label{stsp} If $\eta:A_1 \rar A_2$ is a normal completely
positive map, there exists a standard Hilbert von Neumann $A_1-e_\eta
A_2$ bimodule $\CE_\eta$,  with $e_\eta$ as in Lemma \ref{cp}, which
is singly generated, (i.e., $E = 
[\pi_1(A_1)V\pi_2(e_\eta A_2)]$) with a generator $V \in E$ satisfying
$V^*\pi_1(a_1)V =\pi_2 \circ \eta(a_1)$.

Further, such a pair $(\CE, V)$ of a standard bimodule and generator
is unique in the sense that if $(\widetilde{\CE}, \tilde{V})$ is
another such pair, then there exists $A_i$ - linear unitary operators
$U_i:\CH_i(\eta) \rar \widetilde{\CH_i}, ~i=1,2$ such that $\tilde{V} = U_1
V U_2^*$ and $\widetilde{\CE} = U_1 \CE U_2^*$.  \et

\begin{proof}
  Fix a faithful normal state $\phi$ on $e_\eta A_2$ and set $\CH_2(\eta) =
  L^2(e_\eta A_2,\phi)$, with $\pi_2$ being the left-regular
  representation of $e_\eta A_2$.  We employ the standard notation
  $\hat{a} = \pi(a)\hat{1}$ where $\hat{1}$ is the canonical
  cyclic vector for $\pi(A)$ in $L^2(A)$.  The Hilbert space
  $\CH_1(\eta)$ is obtained after separation and completion of the algebraic
  tensor product $A_1 \otimes e_\eta A_2$ with respect to the semi-inner
  product given by $\langle a_1 \otimes a_2, b_1 \otimes b_2 \rangle =
  \phi(b_2^*\eta(b_1^*a_1)a_2)$; and $\pi_1:A_1 \rar \CL(\CH_1(\eta))$ is
  defined by $\pi_1(a_1)(b_1 \otimes b_2) = a_1b_1 \otimes b_2$. The
  verification that $\pi_1$ is a normal representation is a fairly
  routine application of normality of $\eta$ and $\phi$.

  Define $V: \CH_2(\eta) \rar \CH_1(\eta)$ to be the unique bounded
  operator for which $V(e_\eta \hat{a}_2) = 1 \otimes e_\eta a_2$. 
  For arbitrary $a_1 \in A_1, a_2, b_2 \in e_\eta A_2$, note that
  \beast \langle V^* \pi_1(a_1) V \hat{a}_2, \hat{b}_2 \rangle &=&
  \langle  a_1\otimes a_2, 1 \otimes b_2  \rangle\\\ 
  &=& \phi(b_2^* \eta(a_1) a_2  \rangle\\
  &=& \langle \pi_2(\eta(a_1)) \hat{a}_2,\hat{b}_2 \rangle \eeast thus
  showing that indeed $V^* \pi_1(a_1) V = \pi_2(\eta(a_1))$ for all
  $a_1 \in A_1$.

Set $E = [\pi_1(A_1)V\pi_2(e_\eta A_2)]$ and observe that

\beast [E^*E] &=& [\pi_2(e_\eta
A_2)V^*\pi_1(A_1)\pi_1(A_1)V\pi_2(e_\eta A_2)]\\ 
&=& [\pi_2(e_\eta A_2)\pi_2(\eta(A_1))\pi_2(e_\eta A_2)]\\
&=& [\pi_2(e_\eta A_2 \eta(A_1)e_\eta A_2)]\\
&=& \pi_2(e_\eta A_2) ~,
\eeast
by Lemma \ref{cp}. Further, if $x = \pi_1(a_1)V\pi_2(e_\eta a_2)$ for
$a_i \in A_i$, note that, by definition, we have
$x(\hat{e}_\eta) = a_1 \otimes e_\eta a_2$
and hence, $[\bigcup \{ran~x: x \in E\}] = \CH_1(\eta)$. This shows that
there exist projections $\{p_i:i \in I\} \subset [EE^*]$ such that
$id_{\CH_1(\eta)} = WOT-lim_i p_i$. Hence, we see that
\[\pi_1(A_1) \subset [\bigcup  \{\pi_1(A_1) p_i: i \in I\}] \subset
[\pi_1(A_1)EE^*] \subset [EE^*] ~;\] 
and we have verified everything neeed
to see that the tuple $\CE_\eta = (E, (\pi_1, \CH_1(\eta)), (\pi_2,
\CH_2(\eta)))$ defines a standard Hilbert von Neumann $A_1-e_\eta A_2$ -
bimodule.
   As for the uniqueness assertion, if $(\widetilde{\CE}, \tilde{V})$
  also works, then $\widetilde{\CH_2} =
  L^2(e_\eta A_2,\tilde{\phi})$ for some faithful normal state $\tilde{\phi}$
  on $e_\eta A_2$. In view of the `uniqueness of the standard module of a von
  Neumann algebra' - see [Haa], for instance - there exists an $e_\eta A_2$ - linear unitary
  operator $U_2:\CH_2(\eta) \rar \widetilde{\CH_2}$. Observe next that if
  $\xi, \eta \in \CH_2$ and $a_1, b_1 \in A_1, a_2, b_2 \in e_\eta A_2$, then 
\beast
  \lefteqn{\langle \pi_1(a_1)V\pi_2(a_2)\xi, \pi_1(b_1)V\pi_2(b_2)\eta
    \rangle}\\ 
  &=& \langle \pi_2(b_2^*)V^*\pi_1(b_1^*a_1)V\pi_2(a_2)\xi, \eta \rangle\\
  &=& \langle \pi_2(b_2^*)\pi_2(\eta(b_1^*a_1))\pi_2(a_2)\xi, \eta \rangle\\
  &=& \langle \pi_2(b_2^*\eta(b_1^*a_1)a_2)\xi, \eta \rangle\\
  &=& \langle U_2 \pi_2(b_2^*\eta(b_1^*a_1)a_2)\xi, U_2\eta \rangle\\
  &=& \langle \widetilde{\pi_2}(b_2^*\eta(b_1^*a_1)a_2)U_2\xi, U_2\eta \rangle\\
  &=& \langle  \widetilde{\pi_2}(b_2^*) \widetilde{V^*}
  \widetilde{\pi_1}(b_1^*a_1) \widetilde{V}
  \widetilde{\pi_2}(a_2)U_2\xi, U_2\eta \rangle\\ 
  &=& \langle
  \widetilde{\pi_1}(a_1)\widetilde{V}\widetilde{\pi_2}(a_2)U_2\xi,
  \widetilde{\pi_1}(b_1)\widetilde{V}\widetilde{\pi_2}(b_2)U_2\eta
  \rangle ~.  \eeast Deduce from the above equation and the assumed
  non-degeneracy of $\CE$ and $\widetilde{\CE}$ that there is a unique
  unitary operator $U_1: \CH_1 \rar \widetilde{\CH_1}$ such that
  \be \label{u1def} U_1\left( \pi_1(a_1)V\pi_2(a_2)\xi \right) =
  \widetilde{\pi_1}(a_1)\widetilde{V}\widetilde{\pi_2}(a_2)U_2\xi \ee
  for all $a_1 \in A_1, a_2 \in e_\eta A_2$ and $\xi \in \CH_2(\eta)$
  It is easy to 
  see from equation (\ref{u1def}) that $U_1$ is necessarily $A_1$ -
  linear, that $U_1 V = \tilde{V} U_2$ or $\tilde{V} = U_1 V U_2^*$
  and that $\widetilde{\CE} = U_1 \CE U_2^*$, and the proof of the
  theorem is complete.
\end{proof}

\begin{Remark}\label{eeta1}
Notice that the irritating $e_\eta$ above is equal to the $1$ of $A_2$ in some good cases, such as the following:
\begin{itemize}
\item when $\eta$ is unital, i.e., $\eta(1)=1$;
\item when $\eta(1) \neq 0$ and $A_2$ is a factor.
\end{itemize}
\end{Remark}

The uniqueness assertion in Theorem \ref{stsp} can also be deduced from
the following useful criterion for isomorphism of standard bimodules:

\begin{Lemma}\label{krpisom}
Two standard Hilbert von Neumann $A_2$ bimodules $\CE^{(i)}= (E^{(i)},(\pi^{(i)}_1,\CH^{(i)}_1), (\pi^{(i)}_2,\CH^{(i)}_2)), i=1,2$ are isomorphic if and only if there exist $E_0^{(i)} = \{x_j^{(i)}:j \in I\} \subset E^{(i)}$ such that
\ben
\item $[E_0^{(i)}] = E^{(i)}$, and
\item $(\pi_2^{(1)})^{-1}(x_j^{(1)*}x_k^{(1)}) = (\pi_2^{(2)})^{-1}(x_j^{(2)*}x_k^{(2)}) ~\forall j,k \in I$
\een
\end{Lemma} 

\begin{proof}
The only if implication is clear, as we may choose $E^{(i)}_0 = E^{(i)}$ and $x^{(2)} = U_1x^{(1)}U_2^*$ for all $x^{(1)} \in E^{(1)} (=I)$. Now for the other `if half'.

In view of the `uniqueness of the standard module of a von Neumann algebra - see [Haa] -there exists an $A_2$ - linear unitary operator $U_2:\CH^{(1)}_2\rar \CH^{(2)}_2$. For arbitrary $j,k \in I, \xi_1, \xi_2 \in \CH_2^{(1)}$, observe that
\beast
\langle x_j^{(1)}\xi_1, x_k^{(1)}\xi_2 \rangle &=& \langle \xi_1, x_j^{(1)*}x_k^{(1)} \xi_2 \rangle \\
&=& \langle U_2\xi_1, U_2 \pi_2^{(1)} (\pi_2^{(1)} )^{-1}(x_j^{(1)*}x_k^{(1)}) \xi_2 \rangle \\
&=& \langle U_2\xi_1, \pi_2^{(2)}(\pi_2^{(1)} )^{-1}(x_j^{(1)*}x_k^{(1)})U_2  \xi_2 \rangle \\
&=&  \langle U_2\xi_1, \pi_2^{(2)}(\pi_2^{(2)})^{-1}(x_j^{(2)*}x_k^{(2)})U_2 \xi_2 \rangle \\
&=&  \langle x_j^{(2)}U_2\xi_1, x_k^{(2)}U_2 \xi_2 \rangle ~;
\eeast
deduce from the above equation and the non-degeneracy of the
$\CE^{(i)}$ that there exists a unique unitary operator
$U_1:\CH^{(1)}_1\rar \widetilde{\CH^{(2)}_1}$ such that
$U_1(x_j^{(1)}\xi) =  x_j^{(2)}U_2\xi ~\forall j \in I, \xi \in
\CH_2^{(1)}$. The definitions show that $U_1x_j^{(1)} = x_j^{(2)}U_2
~\forall j \in I$ and hence that $U_1E^{(1)} = E^{(2)}U_2$. Thus
indeed $E^{(2)} = U_1 E^{(1)} U_2^*$ and the proof of the `if half' is complete.

\end{proof}

Notice, incidentally, that in the setting of the  Lemma above, the equation
\[Tx^{(1)} = U_1 x^{(1)} U_2^* \]
defines a WOT-continuous linear bijection $T: E^{(1)}  \rar E^{(2)}$ satisfying
\[Tx^{(1)} (Ty^{(1)})^* Tz^{(1)} = T(x^{(1)} (y^{(1)})^* z^{(1)})\]
for all $x^{(1)}  y^{(1)},  z^{(1)} \in E^{(1)}$.
%\beast
%U_1 (x^{(1)} (y^{(1)})^*) z^{(1)}

\begin{Remark}\label{bimcprmk} {\ben 
\item The `generator' $V$ of Theorem \ref{stsp} is an isometry precisely when $\eta$ is unital.
\item If $\CE$ is a singly
    generated Hilbert von Neumann $A_1-A_2$ bimodule, then it is
    generated by a partial isometry (by Lemma \ref{Epd}). Further, that
    generator, say $V$ may be used to define the obviously completely
    positive map $\eta;A_1 \rar A_2$ by
\[\eta(a_1) = \pi_2^{-1}(V^*\pi_1(a_1)V) ~;\]
and then $\CE$ would be isomosrphic to $\CE_\eta$ if and only if $\CE$
is a standard non-degenerate bimodule. 
\een
}
\end{Remark}

\section{Connes fusion}

\begin{Example}\label{E^d}
{\rm If $\CE = (E,(\pi_1,\CH_1),(\pi_2,\CH_2))$ is a Hilbert von Neumann 
  $A_1-A_2$ - bimodule and $\CK$ is any Hilbert space, then $\CE
\otimes id_\CK = (E \otimes id_\CK,(\pi_1 \otimes 
id_\CK,\CH_1\otimes \CK),(\pi_2 \otimes id_\CK,\CH_2 \otimes \CK))$ is
also a Hilbert von Neumann  $A_1-A_2$ - bimodule, where of course we
write $E \otimes id_\CK$ for $\{x \otimes id_\CK : x \in E\}$.}
\end{Example}

\begin{Lemma}\label{qep}
 Let $\CE = (E,(\pi_1,\CH_1),(\pi_2,\CH_2))$  be a Hilbert von Neumann  
$A_1-A_2$ - bimodule. For  a projection $p \in
\CP(\pi_2(A_2)^\prime)$, let $q$ be the projection with range
$\left[\bigcup\{ran(xp): x \in E\}\right]$. Then 
\ben
\item $q \in \CP(\pi_1(A_1)^\prime)$;
\item $y \in E \Rar qyp = qy = yp$; and
\item $q\CE p = (qEp, (q\pi_1(\cdot), q\CH_1), (p\pi_2(\cdot), p\CH_2))$ satisfies all the requirements for   a non-degenerate Hilbert von Neumann 
  $A_1-A_2$ - bimodule, with the posible exception of injectivity of  $p\pi_2(\cdot)$.
\een

We shall use the suggestive notation $\CE_* p = q$ when $q,\CE,p$ are
so related.
\end{Lemma}

\begin{proof} 1. Since $\pi_1(A_1)E \subset E$, it follows that $ran(q)$
  is stable under $\pi_1(A_1)$.

2. For all $y \in E$, $ran(yp) \subset ran(q) \Rar qyp =
yp$. Next, if $\xi, \eta \in \CH_2$, and $x,y \in E$, note that
\beast \langle xp \xi, y(1-p)\eta \rangle &=& \langle \xi,
px^*y(1-p)\eta \rangle\\  
&\in& \langle \xi, p\left[E^*E\right](1-p)\eta \rangle\\
&=& 0 ~,
\eeast
since $\left[E^*E\right] = \pi_2(A) \subset \{p\}^\prime$; since $\{xp\xi: \xi \in \CH_2\}$ is
total in $ran(q)$, this says that $qy(1-p)=0$, as desired.

3. \be\label{eq1}
\left[(qEp)^*(qEp)\right] = \left[(Ep)^*(Ep)\right] = p\left[E^*E\right]p = p \pi_2(A_2)
\ee
since $\left[E^*E\right] = \pi_2(A) \subset \{p\}^\prime$; while
\be\label{eq2}
\left[(qEp)(qEp)^*\right] = q\left[EE^*\right]q \supset q \pi_1(A_1). \ee

Non-degeneracy of $q\CE p$ follows immediately from equations
(\ref{eq1}) and (\ref{eq2}).
\end{proof}

\begin{Remark}\label{facok}
In general, if $\pi:M \rar \CL(\CH)$ is a faithful normal
representation, and if $p \in \pi(M)^\prime$, the subrepresentation
$p\pi(\cdot)$ is faithful if and only if the central support of $p$ is
$1$ - i.e., $\sup\{upu^*: u \in \pi(M)^\prime\} = 1$.

In particular if the $\CE$ of Lemma \ref{qep} is actually a Hilbert von Neumann 
  $A_1-A_2$ - bimodule, and if $A_2$ happens to be a factor, then the
$q\CE p$ of Lemma \ref{qep} is actually a Hilbert von Neumann
bimodule.
\end{Remark}

We next lead to our description of what is sometimes termed `internal
tensor product' but which we prefer (in view of this terminology being already in use for tensor products of bimodules over von Neumann algebras) to refer to as the {\em Connes
fusion} of Hilbert von Neumann  bimodules. Thus, suppose
$\CE = (E,(\pi_1,\CH_1),(\pi_2,\CH_2))$ is a Hilbert von Neumann 
$A_1-A_2$ - bimodule and $\CF = (F,(\rho_2,\CK_2),(\rho_3,\CK_3))$ is a
Hilbert von Neumann   $A_2-A_3$ - bimodule. We know that the normal
representation $\rho_2$ of $A_2$ is equivalent to a subrepresentation
of an infinite ampliation of the faithful normal representation
$\pi_2$ of $A_2$; thus there exists an $A_2$ - linear isometry $u:\CK_2
\rar \CH_2 \otimes \ell^2$: i.e., $u^*u = id_{\CK_2}$ and $u \rho_2(x) =
(\pi_2(x) \otimes id_{\ell^2}) u ~\forall x \in A_2$. It follows that
$p = uu^* \in (\pi_2(A_2) \otimes id_{\ell^2})^\prime$.   

Now, set $p=uu^*$ and let $q = (\CE \otimes 1_{\ell^2}))_*(p)$ be
associated to this $p$ as in Lemma \ref{qep} (applied to $\CE \otimes
1_{\ell^2}$). 

Finally, if $x \in E, y \in F$, define $x \bigodot y$ to be the
composite operator
\[\CK_3 \stackrel{x \bigodot y}{\longrightarrow} q(\CH_1 \otimes
\ell^2) ~=~ \CK_3 \stackrel{y}{\longrightarrow} \CK_2
\stackrel{u}{\longrightarrow} uu^*(\CH_2 
\otimes \ell^2) \stackrel{x \otimes id_{\ell^2}}{\longrightarrow} q(\CH_1
\otimes \ell^2) ~,\]
set $E \bigodot F = \left[\{x \bigodot y:x \in E, y \in F\}\right]$; and
finally define the Connes fusion of $\CE$ and $\CF$ to be
\be \label{cfusdef}
\CE \otimes_{A_2} \CF = (E \bigodot F, (q(\pi_1 \otimes
id_{\ell^2})|_{ran~q},q(\CH_1 \otimes \ell^2)), (\rho_3,\CK_3))~. 
\ee

The justification for our use of `Connes fusion'  for our construction
lies (at least for  standard bimodules, by Lemma \ref{krpisom}) in the fact  that (in the notation
defining Connes fusion) the $A_3$ - valued inner 
product on $\CE \circ \CF$ satisfies 
\beast 
\langle x_1 \bigodot y_1, x_2 \bigodot y_2 \rangle_{A_3}
&=& (x_1 \bigodot y_1)^*(x_2 \bigodot y_2)\\
&=& (x_1 \otimes id_{\ell^2})uy_1)^*(x_2 \otimes id_{\ell^2})uy_2\\
&=& y_1^*u^*(x_1^*x_2 \otimes id_{\ell^2})uy_2\\
&=& y_1^*(x_1^*x_2)y_2 ~~~\mbox{ (since $u$ is an $A_2$ - linear isometry)}\\
&=& y^*_1 \langle x_1, x_2 \rangle_{A_2} y_2\\
&=& \langle y_1, \langle x_1, x_2 \rangle_{A_2} y_2 \rangle_{A_3}~.
\eeast

\begin{Proposition}\label{cfusprop}
The Connes fusion of (non-degenerate) Hilbert von Neumann  bimodules is
again a (non-degenerate) Hilbert von Neumann  bimodule.
\end{Proposition}

\begin{proof} 
Clearly $E \bigodot F$ is a WOT-closed linear space of operators between the
  asserted spaces. Observe next that

\beast 
\left[(E \bigodot F)(E \bigodot F)^*\right] &=& \left[\{((x_1 \otimes
  id_{\ell^2})uy_1)(x_2 \otimes id_{\ell^2})uy_2)^*: x_i \in E, y_j \in F\}\right]\\
&=& \left[\{(x_1 \otimes id_{\ell^2})uy_1y_2^*u^*(x_2 \otimes
      id_{\ell^2}): x_i \in E, y_j \in F\}\right]\\
&=& \left[\{(x_1 \otimes id_{\ell^2})u\left[FF^*\right]u^*(x_2 \otimes id_{\ell^2}):
    x_i \in E\}\right]\\
&\supset& \left[\{(x_1 \otimes id_{\ell^2})u \rho_2(A_2) u^*(x_2 \otimes id_{\ell^2}):
    x_i \in E\}\right]\\
&=& \left[\{(x_1 \otimes id_{\ell^2})(\pi_2(A_2) \otimes
    id_{\ell^2})uu^* (x_2 \otimes id_{\ell^2}): x_i \in E\}\right]\\
&=& \left[(E \otimes id_{\ell^2}) uu^* (E \otimes id_{\ell^2})^*\right]
  ~~~(\mbox{since } E \pi_2(A_2) = E)\\
&=& q(\pi_1(A_1) \otimes id_{\ell^2})
\eeast
(in particular $q \in \left[(E \bigodot F)(E \bigodot F)^*\right]$) and that
\beast
\left[(E \bigodot F)^*(E \bigodot F) \right] &=& \left[\{((x_1 \otimes
  id_{\ell^2})uy_1)^*(x_2 \otimes id_{\ell^2})uy_2): x_i \in E, y_j \in F\}\right]\\
&=& \left[\{(y_1^*u^*(x_1^*x_2 \otimes id_{\ell^2}))uy_2): x_i \in E, y_j \in F\}\right]\\ 
&=& \left[\{(y_1^*u^*(\pi_2(A_2) \otimes id_{\ell^2})uy_2): y_j \in F\}\right]\\ 
&=& \left[\{(y_1^*u^*u\left[\rho_2(A_2)\right]y_2): y_j \in F\}\right]\\ 
&=& \left[\{(y_1^*(\rho_2(A_2))y_2): y_j \in F\}\right]\\ 
&=& F^*F~~~(*)\\ 
&=& \rho_3(A_3)~,
\eeast
where the justification for the step labelled (*) is that $\rho_2(A_2)F
= F$ (see Remark \ref{HvNmodrmks} (7). This completes the verification
that $\CE \otimes_{A_2}\CF$ is indeed a Hilbert von Neumann $A_1-A_3$ bimodule.

Now, suppose $\CE$ and $\CF$ are both non-degenerate. Then
\beast
\lefteqn{\xi \in \bigcap\{ker z : z~ in~ E \bigodot F\}}\\
& \Rar &  (x \otimes id_{\ell^2})uy \xi = 0 ~\forall x \in E, y \in F\\
& \Rar &  uy \xi = 0 ~\forall y \in F ~~~\mbox{(as $E \otimes
  id_{\ell^2}$ is non-degenerate)}\\
& \Rar &  y \xi = 0 ~\forall y \in F ~~~\mbox{(as $u$ is isometric)}\\
& \Rar &  \xi = 0 ~~~\mbox{(as $F$ is non-degenerate)} ~;
\eeast
while
\beast
\lefteqn{\left[\bigcup\{ran((x \otimes id_{\ell^2})uy) : x \in E, y \in F\}\right]}\\
&=& \left[\bigcup\{ran((x \otimes id_{\ell^2})u) : x \in E\}\right] ~~~(\mbox{since
  $F$ is non-degenerate)}\\
&=& \left[\bigcup\{ran((x \otimes id_{\ell^2})uu^*) : x \in E\}\right]\\
&=& ran~q ~~~\mbox{(by ~definition)}
\eeast
and hence $E \bigodot F$ is indeed non-degenerate.

\end{proof} 

Before addressing the question of the dependence of the definition of
Connes fusion and the seemingly {\em ad hoc} $A_2$ - linear partial 
isometry $u$, we introduce a necessary definition and the ubiquitous lemma.

\begin{Definition} \label{bimisom} Two Hilbert von Neumann $A_2$ modules,
  say $\CE^{(i)} = (E^{(i)}, \CH_1^{(i)}, (\pi^{(i)}_2,\CH_2^{(i)})), i=1,2$ 
  are considered {\em isomorphic} if there exists unitary
  operators $w_j: \CH^{(1)}_j \rar \CH^{(2)}_j$, with $w_2$ being $A_2$ - linear, such that
\[E^{(2)} = w_1 E^{(1)} w_2^* ~.\]

If the $\CE^{(i)}$ happen to be $A_1-A_2$ bimodules, they are said to be isomorphic if, in addition to the above, the unitary $w_1$ happens to be $A_1$ - linear.
\end{Definition}

\begin{Lemma}\label{qepw}
Let $\CE = (E,(\pi_1,\CH_1),(\pi_2,\CH_2))$ be a Hilbert von Neumann $A_1-A_2$ bimodule.
Suppose $w \in \pi_2(A_2)^\prime$ is a partial isometry with $w^*w =
p, ww^* = \tilde{p}$. Let $q = \CE_*p$ and $\tilde{q}=\CE_*\tilde{p}$ in the notation of
Lemma \ref{qep}. Then there exists a unique partial isometry $w_1 \in
\pi_1(A_1)^\prime$ such that $w_1^*w_1 = q, w_1w_1^* = \tilde{q}$. 
\end{Lemma}

\begin{proof}
We first assert that there is a unique unitary operator $W_1:q(\CH_1)
\rar \tilde{q}(\CH_1)$ satisfying $WTp = Tw~ \forall T \in E$. This is
because:
\begin{itemize}
\item $(T_1w)^*(T_2w) = w^*T_1^*T_2w = T_1^*T_2 p = p^*T_1^*T_2p,
  ~\forall T_1, T_2 \in E$ and
\item $q(\CH_1) = \left[\bigcup \{ran(Tp):T \in E\}\right]$ and
$\tilde{q}(\CH_1) = \left[\bigcup \{ran(Tw):T \in E\}\right]$ (since $ran~w
  = ran~\tilde{p}$.
\end{itemize}
Finally $w_1 = W_1q$ does the job.
\end{proof}

\begin{Remark}\label{cfwdef}
\ben
{\rm \item
We now verify that the definition we gave of $\CE \otimes_{A_2}\CF$ is really independent of the choice of the isometry $u$ used in that definition. Indeed, suppose $u, \tilde{u}: \CK_2 \rar \CH_2 \otimes \ell^2$ are two $A_2$ - linear isometries. If $uu^* = p, \tilde{u}\tilde{u}^* = \tilde{p}$, then $w = \tilde{u}u^*$ is a partial isometry in $(\pi_2(A_2) \otimes id_{\ell^2})^\prime$ with $w^*w = p, ww^* = \tilde{p}$.
Now apply Lemma \ref{qepw} to $\CE \otimes id_{\ell^2}$ and $w,p,\tilde{p}$ to find a $W \in (\pi_1(A_1) \otimes id_{\ell^2})^\prime$ such that $W^*W = q = (\CE \otimes id_{\ell^2})_*p$ and
$WW^* = \tilde{q}  = (\CE \otimes id_{\ell^2})_*\tilde{p}$. Then, as the proof of Lemma \ref{qepw} shows,
$W:q(\CH_1 \otimes \ell^2) \rar \tilde{q}(\CH_1 \otimes \ell^2)$ is a unitary operator satisfying $W(x\otimes id_{\ell^2})p =( x \otimes id_{\ell^2}) w~ \forall x \in E$. It is now a routine matter to verify that the unitary operators $W: q(\CH_1 \otimes \ell^2) \rar \tilde{q}(\CH_1 \otimes \ell^2$ and $id_{\CK_3}$ establish an isomorphism between the models of  $\CE \otimes_{A_2} \CF$ given by $u$ and $\tilde{u}$ are isomorphic.
\item A not dissimilar reasoning shows that the isomorphism type of the Connes fusion of teo standard bimodules depends only on the isomorphism classes of the two `factors' in the fusion, and is also standard.
\item If $\CE$ is only a Hilbert von Neumann $A_2$-module, and $\CF$
  is a Hilbert von Neumann $A_2-A_3$-bimodule, their Connes fusion
  $\CE \otimes_{A_2} \CF$ would still make sense as a Hilbert von
  Neumann $A_3$-module.
}
\een
\end{Remark}

\section{Examples}

We now discuss some examples of Hilbert von Neumann (bi)modules. 
\ben
\item The simplest (non-degenerate) example is obtained when $A_j =
  \CL(\CH_j), \pi_j = 
  id_{A_j}$ for $j=1,2$ and $E = \CL(\CH_2,\CH_1)$; all the
  verifications reduce just to matrix multiplication.
\item Suppose $A_2$ is a unital von Neumann subalgebra of $A_1$, and
  suppose there exists a faithful normal conditional expectation
  $\e:A_1 \rar A_2$. Let $\phi_2$ be a faithful normal state (even
  semi-finite weight will do). Let $\phi_1 = \phi_2 \circ \e, 
  \CH_j = L^2(A_j,\phi_j)$, and let $\pi_j$ be the left regular
  representation of $A_j$ on $\CH_j$. Write $U$ for the natural
  isometric identification of $\CH_2$ as a subspace of $\CH_1$ (so
  that the `Jones projection' will be just $UU^*$). Finally, define
\[\CE_{(A_2 \subset A_1)} = (\pi_1(A_1)U,
(\pi_1,\CH_1),(\pi_2,\CH_2))\]
In this case, we find that $[EE^*] = [\pi_1(A_1)e\pi_1(A_1)]$, and we
find the `basic construction of Jones appearing naturally in this
context.

Further, it is a consequence of the uniqueness assertion in Theorem
\ref{stsp} that $\CE_\e \cong \CE_{(A_2 \subset A_1)}$.

\item Suppose $(M, \CH, J, P)$ is a standard form of $M$ in the sense of [Haa]. 
As indicated in [Haa], there is a canonical `implementing' unitary representation
\[Aut(M) \ni \theta \mapsto u_\theta \in \CL(\CH)\]
satisfying $u_\theta x u_\theta^* = \theta (x) ~\forall x \in M$. We have the natural Hilbert von Neumann $M-M$ bimodule given by
\[\CE_\theta = (M u_\theta, (id_M, L^2(M)), (id_M,L^2(M)))\]
\item If $\theta,\phi \in Aut(M), M$ are as in the previous example, we
  see now that `Connes fusion corrsponds to composition' in this
  case:
\[\CE_\theta \otimes_M \CE_\phi \cong \CE_{\theta\phi}\]
({\em Reason:} The `$u$' in the definition of Connes fusion is just $id_M$, while 
\[Mu_\theta Mu_\phi = M\theta(M)u_\theta u_\phi = M u_{\theta\phi}~.)\]
\een

\begin{Proposition}
If $\theta,\phi \in Aut(M)$ are as in Example (4) above, then $\CE_\theta \cong \CE_\phi$ if and only if $\theta$ and $\phi$ are inner conjugate.
\end{Proposition}
\begin{proof} 
First, note that any $M$-linear unitary operator on $L^2(M)$ has the form $Jv^*J$ for some unitary $v \in M$, where of course $J$ denotes the modular conjugation operator. Observe next that each $u_\theta$ commutes with $J$ since $\theta$ is a *-preserving map, and hence, for any $x \in M$, we have 
\be \label{uthcom} u_\theta Jv^*J  = J\theta(v^*)J u_\theta  \ee

If $\CE_\theta$ is isomorphic to $\CE_\phi$, there must exist unitary $v_1,v_2 \in M$ such that
\beast Mu_\phi &=& Jv_1^*J Mu_\theta Jv_2J\\
&=& M Jv_1^*J u_\theta Jv_2J\\
&=& M Jv_1^*J  J\theta(v_2)J u_\theta\\
&=& M Jv_1^*\theta(v_2)J u_\theta~;
\eeast
in particular, there must exist a $y \in M$ such that
\[u_\phi = y Jv_1^*\theta(v_2)J u_\theta~.\]
We find that $y$ is necessarily unitary
and hence, writing $u$ for $y$ and $v$ for $ v_1^*\theta (v_2)$, we see that there must be a unitary
$u \in M$ such that
\beast \phi(x) &=& u_\phi x u_\phi^*\\
&=& u JvJ u_\theta x u_\theta^* Jv^*J u^*\\
&=& u JvJ \theta (x)  Jv^*J u^*\\
&=& u  \theta (x)  u^*~.
\eeast
In other words, $\phi$ and $\theta$ are indeed inner conjugate.

Conversely, if $\phi(\cdot) = u \theta (\cdot) u^*$ for some unitary $u \in M$, we see that
$u_\phi = u JuJ u_\theta = u u_\theta J \theta^{-1}(u)J$; so we find that $w_1 = id_M$ and $w_2 = J\theta^{-1}(u)^*J$ define $M$ - linear unitary operators on $L^2(M)$ such that
$M u_\phi = Mu u_\theta w_2^* = w_1 Mu_\theta w_2^*$, thereby
establishing that $\CE_\theta \cong \CE_\phi$. 
\end{proof}

\bigskip \noindent
{\bf {\large References}}

\medskip \noindent
[Haa] Haagerup, U. {\em The standard form of von Neumann algebras}, Math. Scand., {\bf 37}, 1975, 271-283.

\medskip \noindent
[Lan] Lance, Christopher, {\em Hilbert $C^*$ - algebras}, LMS Lecture
Note Series 210, CUP, Cambridge.

\medskip \noindent
[Skei] Skeide, Michael, {\em Hilbert modules and applications in                
  quantum probability}, Cottbus, 2001.

\end{document}